\magnification=1200
\documentstyle{amsppt}

\TagsOnRight

\define\C{{\Bbb C}}
\define\R{{\Bbb R}}
\define\re{\text{\rm Re}\,}
\define\im{\text{\rm Im}\,}
\define\aut{\text{\rm Aut}\,}
\define\st{\,|\,}
\define\E{{\Cal E}}

\topmatter

\title
Factorization of proper holomorphic \\
mappings through Thullen domains
\endtitle

\author
Kang-Tae Kim, Mario Landucci, Andrea F. Spiro
\endauthor

\address
(Kim) Department of Mathematics,
Pohang University of Science and Technology, Pohang 790-784
Republic of Korea (South) 
\endaddress
\email
kimkt\@postech.ac.kr
\endemail

\address
(Landucci and Spiro) Dipartimento di Matematica,
Universit\`a di Ancona, Ancona, Italy
\endaddress
\email
landucci\@anvax1.cineca.it, spiro\@anvax1.cineca.it
\endemail

\thanks
Research of the first named author is supported in part by
Grant POSTECH/BSRI-95-1429, POSTECH/BSRI Special Fund, GARC
and KOSEF Grant 961-0104-021-2 of the Republic of Korea (South).  
Research of the second and third named authors is supported by
M.U.R.S.T. 40\% and 60\% funds.
\endthanks

\leftheadtext{K.T. Kim, M. Landucci and A.F. Spiro}
\rightheadtext{Factorization of Proper Holomorphic Maps}

\abstract
Given a bounded domain in $\C^2$ and a proper holomorphic mapping
onto the unit ball in $\C^2$, we give a criterion for a 
proper mapping to factor through Thullen domains by 
its branch locus.  As an application, we present 
a characterization of the Thullen domains among bounded
domains admitting a proper holomorphic mapping onto the unit ball.
\endabstract


\endtopmatter

\document

\head
1. Introduction
\endhead

A continuous map is called \it proper\rm, if the pre-image of every 
compact set is compact.  Let $f:D \to G$ be a proper holomorphic
mapping between bounded pseudoconvex domains in $\C^n$.  Denote by 
$Z_{df}$ the analytic set defined by the relation that its holomorphic
Jacobian determinant is zero.  The Proper Mapping 
Theorem of Remmert and Stein (\cite{10}) implies that $f$ is
a finite-to-one holomorphic covering map from $D \setminus Z_{df}$
onto its image.  Thus, all proper holomorphic mappings between 
equi\-dimensional bounded domains of holomorphy possess a generic degree.  

For every positive integer $k$, we consider the Thullen domain
$$
\E_k = \{(z,w) \in \C^2 \st |z|^{2k} + |w|^2 < 1\}.
$$
For the entire paper, we use the notation
$$
\varphi_k (z,w) = (z^k, w) : \E_k \to B^2
$$
where $B^2$ is the open unit ball in $\C^2$.

The primary goal of this article is to establish
the factorization theorem presented in the following.  This theorem
shows in particular that the Thullen domains $\E_k$ together 
with the map $\varphi_k$ form a standard model for proper mappings 
from a bounded domain with a real analytic boundary to the unit ball 
in $\C^2$, in the sense that a proper mappings with controlled
branch locus factors through the map $\varphi_k$ for an appropriate 
positive integer $k$.

\proclaim{Theorem 1} \it 
Let $D$ be a bounded simply connected pseudoconvex domain in $\C^2$ 
with a real analytic boundary.  Assume that $f:D \to B^2$ is a 
proper holomorphic mapping with generic degree $m$ 
such that the analytic variety $Z_{df}$ admits an 
irreducible component $V$ satisfying: 
\roster
\item $f^{-1}(f(V \cap \partial D)) = V \cap \partial D$.
\item $V \cap \partial D$ is connected and contains no singular
point of the variety $Z_{df}$.
\endroster
Denote by $k$ the degree of the mapping $f$ in a sufficiently
small tubular neighborhood of $V$.  Then $m$ is divisible by $k$ and there 
exists a generically $m/k$-to-1 proper holomorphic mapping 
$g:D \to \E_k$ which extends holomorphically 
across $\partial D$ such that
$$
f = \alpha \circ \varphi_k \circ g
$$
for some holomorphic automorphism $\alpha$ of $B^2$, where 
$\varphi_k (z,w) = (z^k,w)$.
\endproclaim

By [9] (see also [2], [3], [7]) all proper holomorphic 
mappings between bounded pseudoconvex domains with entirely real 
analytic smooth boundaries extend holomorphically across the boundaries. 
Thus, the conditions (1) and (2) in the statement of Theorem 1 
(as well as in the statement of Theorem 2 below) are meaningful.  

We would like to point out that Theorem 1 treats the factorization
of proper holomorphic mappings \it into \rm the ball.  This distinguishes
our theorem from a result by Rudin ([11]) which study
the factorization of proper maps \it from \rm the ball.

Notice that Theorem 1 implies in particular the following characterization 
of the Thullen domain $\E_k$ among the source domains of proper holomorphic 
mappings onto the unit ball by its branch locus.

\proclaim{Theorem 2} 
\it Let $D$ be a bounded simply connected pseudoconvex domain in $\C^2$
with a real analytic boundary.  If $D$ admits a generically
$k$-to-1 proper holomorphic mapping $f:D \to B^2$ such that there  
exists an irreducible subvariety $V$ of $Z_{df}$ satisfying:
\roster
\item $f$ is a local $k$-to-1 branched covering at every point of 
$V \cap \partial D$;
\item $V \cap \partial D$ is connected and contains no singular
point of the variety $Z_{df}$;
\endroster
then $D$ is biholomorphic to $\E_k$.
\endproclaim

We would like to remark that the hypothesis in the factorization theorem
above is essential in the sense that there are many 
examples of unfactorizable proper mappings with branch
locus violating the conditions in the hypothesis.  

Although it is desirable in principle to obtain a generalization of the
above theorems to higher complex dimensions, one of the main difficulties 
lies in the fact that there are no known effective normal
forms of Chern-Moser type near the weakly pseudoconvex points in complex 
dimensions higher than two.  Since our methods depend heavily upon the
Chern-Moser normal forms and subsequent developments and generalizations
by Barletta and Bedford ([4]), our two dimensional theorem seems in fact 
optimal.  

\remark{Acknowledgments}  
We are indebted to the referee for valuable suggestions especially
concerning 4.3 Theorem.
\endremark

\head
2. Organization of the paper
\endhead

The rest of this article is organized as follows. In order to 
have the exposition as clear as possible, Section 3 presents the
precise definitions and standard notation for the concepts such 
as pointed CR surfaces 
and associated domains, holomorphic mappings between them, branch
points, local automorphisms and weakly spherical normal forms of
Barletta and Bedford.  Section 4 discusses the rigidity of holomorphic
mappings from a normalized weakly spherical surface into a normalized
strongly pseudoconvex surface in the sense of Chern-Moser.  The first 
key step (4.1 Proposition) is an involved investigation on the effect
of holomorphic mappings to the normal forms.  Based upon this result,
we then show that the only weakly spherical pointed CR surface that 
admits a regular branched holomorphic mapping to a pointed Siegel surface
is defined by $\re w - |z|^{2k} = 0$.  We further show that the 
holomorphic mapping
has to be essentially $(z,w) \mapsto (z^k,w)$.  This establishes
an important step towards our proof of the factorization theorem.  In
Section 5, the proofs of main theorems are given.  The arguments
rely upon the local rigidity established in Section 4 and
an analysis of the splitting of holomorphic correspondences.

\head
3. Basic terminology and notation
\endhead

\subhead
3.1 Pointed CR hypersurfaces and associated domains
\endsubhead
Let $(M,p)$ be a real three dimensional real analytic pointed 
hypersurface in $\C^2$. It naturally admits the standard CR
structure.  Furthermore, the Implicit Function Theorem shows that 
there
exists a real-valued real analytic function $\rho_M (z,w)$ defined
in an open neighborhood $U_M$ of $p$ in $\C^2$ such that 
$$
M \cap U_M = \{(z,w) \in U_M \st \rho_M (z,w) = 0\}    
$$
and such that the gradient $\nabla \rho_M$ is not zero at any point of
$M \cap U_M$.  Shrinking $U_M$ sufficiently small, we may assume that
$U_M \setminus M$ consists of two domains.  We will call such $U_M$
a \it neighborhood associated with \rm $(M,p)$.  In case $M$ is 
pseudoconvex and not Levi flat, only one of the two is
a pseudoconvex domain. We assume that
$$
\Omega_M \equiv \{(z,w) \in U_M \st \rho_M (z,w) > 0\}
$$
is pseudoconvex and call it a \it domain associated with $(M,p)$. \rm
Although the concepts of associated neighborhoods and
associated domains contain ambiguities, they can be readily made
precise by exploiting the concept of germs.  Therefore, in this paper
we choose to use the above terminologies without introducing the 
concept of germs, since there is no danger of confusion in this
exposition.

\subhead
3.2 Holomorphic mappings of pointed CR hypersurfaces
\endsubhead
The next concept to consider is the mapping between two pointed 
CR hypersurfaces (or equivalently, between their germs) in $\C^2$.
Given two pointed CR hypersurfaces $(M,p)$, $(N,q)$ with 
associated neighborhoods $U_M$, $U_N$ and associated domains 
$\Omega_M$,
$\Omega_N$ respectively, we say a mapping $F:(M,p) \to (N,p)$ to be
\it a holomorphic mapping between pointed CR hypersurfaces\rm, if it
extends to a holomorphic mapping $F:U_M \to U_N$ satisfying the 
conditions
$$
F(p)=q, \quad F(M)\subset N, \quad F(\Omega_M) \subset \Omega_N.
$$

We also introduce the notation for the branch locus (or
equivalently, the Segre variety) of $F$ in the following:
$$
Z_{dF} \equiv \{z \in U_M \st \det JF(z) = 0 \}. 
$$
By a \it branch point \rm of the mapping $F$ we mean an 
element of $Z_{dF}$.  We call a branch point \it regular\rm, if
it is a smooth point of $Z_{df}$.

\subhead
3.3 Local Automorphisms
\endsubhead
Furthermore, for a pointed CR hypersurface $(M,p)$ we denote by
$$
\aut (M,p) = \{ \varphi:(M,p) \to (M,p) \st \hbox{1-1 and 
holomorphic}\}
$$
the \it local automorphism group for $M$ at $p$. \rm

\subhead
3.4 Weakly spherical hypersurfaces and normal forms
\endsubhead
Finally, we introduce an important special class of pointed CR 
hypersurfaces following \cite{4}.  A real analytic CR hypersurface 
$M$ is called {\it weakly spherical of order $k_0 \ge 1$ at $p \in 
M$}, if it admits a local holomorphic coordinate system 
$(z,w)$ at $p$ ($p=(0,0)$ in the new coordinates) in which 
$(M,p)$ is defined by the following normalized defining function
(an analogue of the Chern-Moser normal form introduced in [8])
$$
\rho(z,w) = \re w - |z|^{2k_0} - 
\re \sum_{j,h\ge 1 \atop j+h\ge 2k_o +1} 
   a_{jh} (\im w) \, z^j \bar z^h                 \tag3.4.1
$$
where $a_{jh}$ are real analytic functions in the variable
$\im w$ satisfying:
$$
\align
a_{2k_0\, 2k_0} & = a_{3k_0\, 3k_0} = 0 \\
a_{jk_0} = a_{k_0 j} & = 0, \forall j \ge k_0 + 1   \tag3.4.2
\endalign
$$

If a real analytic CR hypersurface $M$ in $\C^2$ contains the 
origin $o=(0,0)$, and if a defining function for $M$ in a 
neighborhood of $o$ takes the form (3.4.1) without a change of
coordinates of $\C^2$ and satisfies the
condition (3.4.2), then we call the pointed CR hypersurface $(M,o)$ 
\it normalized weakly spherical.  \rm  

For a pointed real analytic CR hypersurface $(M,o)$, we say that a 
defining function is \it reduced \rm if it is of the form 
$\rho(z,w) = \re w - G(z,\bar z, \im w)$.  Note that 
in each fixed holomorphic coordinate system, a CR hypersurface
admits at most one reduced defining function.

Suppose that there exists a holomorphic mapping $F$ from $(M,o)$ into
a smooth strongly pseudoconvex CR hypersurface $(N,o)$ for which
$o$ is a regular branch point.  Then, according to [4], for every
\it reduced \rm defining function $\rho_M$ for $(M,o)$, there exists
a local biholomorphism $\gamma$ at the point $o$ 
such that $\rho_M\circ \gamma$ is the weakly spherical normal form, 
which takes the form (3.4.1) and satisfies (3.4.2). See [4] for 
detailed arguments.

\head
4. Holomorphic mapping into the Siegel hypersurface
\endhead
We now investigate the actions of the holomorphic mappings on the 
normal forms of the CR hypersurfaces.  We first present the following
technical result on the holomorphic mappings from a normalized 
weakly spherical CR hypersurface to a normalized strongly pseudoconvex 
hypersurface.

\proclaim{4.1 Proposition}
\it
Let $(M,o)$ be a normalized pointed weakly spherical CR hypersurface
in $\C^2$ of order $k_0 > 2$, and let $(N,o)$ be a normalized
CR hypersurface in $\C^2$ which is strongly pseudoconvex at the
origin $o=(0,0)$.  Then, for every holomorphic mapping
$F:(M,o) \to (N,o)$ with a regular branching point at $o$,
there exist local automorphisms $g_N \in \aut (N,o)$
and $g_M \in \aut(M,o)$ 
 such that
$$
\rho_M = \rho_N \circ g_N \circ F \circ g_M
$$
where $\rho_M$ and $\rho_N$ are the defining functions in
normal form for $(M,o)$ and $(N,o)$, respectively.
\endproclaim


\remark{\bf 4.2 \it Remark} 
A direct application of Hopf's lemma implies only that
$\rho_N \circ F$ is also a defining function of $(M,o)$, which in 
turn 
yields at best that $h\cdot\rho_M = \rho_N \circ F$ for some positive
real valued function $h$.  The first merit of the proposition
above is that we have limited the ambiguity to a significant degree, 
by showing that $h$ can be replaced by a composition by two local 
automorphisms of the  surfaces.   Even with this fact alone,
it seems providing an interesting aspect for a further investigation, 
which should be of a separate interest.
\endremark

\demo{Proof of 4.1 Proposition}
Let $W_M$ be the set of weakly pseudoconvex points on $M \cap U_M$.
It is given by 
$$
W_M = \{(z,w) \in M \cap U_M \st z = 0\}
$$
since $M$ is normalized.  Since $F$ maps holomorphically $M$ into $N$
and $\Omega_M$ into $\Omega_N$ respectively, Hopf's lemma implies that 
$$
W_M = \{(z,w) \in M \cap U_M \st \det(JF)(z,w) = 0 \}.
$$
Since $o$ is a regular branch point for $F$, the Weierstrass
Preparation Theorem implies that, by shrinking $U_M$ if necessary,
the mapping $F:U_M \to F(U_M) \subset U_N$ is a branched covering 
such that $F(Z_{df} \cap U_M)$ is a smooth analytic variety.  
Furthermore, we may choose a change of holomorphic local coordinates 
at $o$ by a local biholomorphism $\psi$ so that
$$
\align
(\psi \circ F)(Z_{df} \cap U_M) 
          & \subset \{(0,w) \in \C^2 \st w \in \C\}   \\
(\psi \circ F)(0,w) & = (0,w) 
\endalign
$$

Hence, 
$$
(\psi \circ F)(z,w) = (z^m\cdot k(z,w), w+z\cdot g(z,w))
$$
where $k(0,w)$ is not identically zero.  Since the map
$(z,w) \mapsto (z,w+z \cdot g(z,w))$ is invertible near the origin
$o=(0,0)$, we define a local biholomorphism $\phi_1$ in a neighborhood
of the origin by
$$
\phi_1^{-1} (z,w) = (z, w+ z \cdot g(z,w)).
$$
Thus we have
$$
(\psi \circ F \circ \phi_1)  (z,w) = (z^m h(z,w), w),
$$
where $h(z,w) = (k\circ\phi_1)(z,w)$.
Now we show that $h(0,0) \not= 0$. Since
$$
\det J(\psi \circ F \circ \phi_1) 
= z^{m-1}\left( mh + z \frac{\partial h}{\partial z} \right)   \tag4.1
$$
and since the origin $o$ is a smooth point of the variety
$Z_{d(\psi\circ F\circ \phi_1)}$, the equation 
$\det J(\psi \circ F\circ\phi_1) = 0$ must determine a smooth 
variety near the origin.  From (4.1) above, we must then have
either $h(0,0) \not= 0$ or $mh + z \frac{\partial h}{\partial z}$
vanish identically on $\{(0,w) \in \C^2\}$.  However the latter implies
that $h(0,w)$ is identically zero.  Since $h(0,w)=k(0,w)$, this 
contradicts to the fact that $k(0,w)$ is not identically zero.
Hence $h(0,0) \not= 0$.  (This paragraph is a modification of 
an argument in Barletta-Bedford \cite{4}.)

Therefore, in an open neighborhood of the origin we may choose a
local biholomorphism $\phi_2$ given by
$$
\phi_2^{-1} (z,w) = (z (h(z,w))^{1/m}, w)
$$
by taking any branch.  Let 
$$
\phi = \phi_1 \circ \phi_2 .
$$
Then we have  
$$
\tilde F (z,w) \equiv (\psi \circ F \circ \phi)(z,w) = (z^m, w)     
\tag4.2
$$
and
$$
F^\#(z,w) \equiv (\psi \circ F)(z,w) = (\tilde F \circ 
\phi^{-1})(z,w)
=
(z^m k(z,w), w + z\cdot g(z,w)) 
\tag4.3
$$
in a neighborhood of the origin $o$. \par

We now consider the pointed CR hypersurfaces $(N',o)$ and $(M',o)$ 
defined by 
$$
N'\equiv \psi(N\cap U_N),  
\qquad 
M' \equiv \phi^{-1}(M\cap U_M) 
= \tilde F^{-1}(N'\cap \psi(U_{N})),
$$
where $U_M$ and $U_{N}$ are taken sufficiently small so that 
$\psi$ and $\phi$ are well-defined.

To conclude the proof, we prove the following 

\proclaim{Claim}
\it  There exist two
\underbar{reduced} defining functions $\rho_{N'}$ and $\rho_{M'}$ for 
$(N',o)$ and $(M',o)$ such that 
$$
\rho_{M'} = \rho_{N'}\circ \tilde F
$$
\endproclaim

Suppose for the moment that the claim is proved.  Since $\rho_N$
and $\rho_M$ are two defining functions in normal
form, according to \cite{4}  and \cite{8} there exist 
two local biholomorphic mappings $\gamma: (N',o) \to (N,o)$ and 
$\chi:(M',o) \to (M,o)$ such that 
$$
\rho_N = \rho_{N'}\circ \gamma^{-1}, \qquad 
\rho_M = \rho_{M'}\circ \chi^{-1}
$$
This implies that
$$
\align
\rho_M & = \rho_N \circ \gamma \circ \tilde F \circ \chi^{-1} \\
           & = \rho_N \circ (\gamma\circ \psi) 
                   \circ F \circ (\phi\circ \chi^{-1}),
\endalign
$$
where it is obvious that $(\gamma\circ \psi) \in \aut(N,o)$ and
$(\phi\circ\chi^{-1}) \in \aut(M,o)$.  Thus the assertion of
the proposition is obtained.  \par

Therefore it remains only to prove the claim above.
Let us now look for a \it reduced \rm defining function $\rho_{N'}$.
If we denote by
$$
\delta = \rho_N \circ \psi^{-1}
$$
then, by shrinking $U_M$ if necessary, we have
$$
F^\# (\Omega_M) = \{(z,w) \in U_M \st \delta (z,w) < 0\}  \tag4.4
$$
and
$$
F^\# (M \cap U_M) = \{(z,w) \in U_M \st \delta (z,w) = 0\}. \tag4.5
$$

By Hopf's Lemma, $\delta \circ F^\#$ is also a defining function 
for $(M,o,\Omega_M)$.  Therefore,
$$
\delta \circ F^\# (z,w) = \alpha(z,w) \, \rho_M (z,w)
$$
for some positive real analytic function $\alpha:U_M \to \R$.

Write as $z = x+ iy$ and $w = u+iv$.  Since $F^\#(z,w) = 
(z^m k(z, w), w + z\cdot g(z,w))$,
we have
$$
\left. \frac{\partial(\delta\circ F^\#)}{\partial u} \right|_o
= \left.\frac{\partial \delta}{\partial u} \right|_o .  \tag4.6
$$
On the other hand,
$$
\left.\frac{\partial(\delta\circ F^\#)}{\partial u} \right|_o
= \left.\frac{\partial (\alpha\cdot\rho_M)}{\partial u}\right|_o
= \left.\alpha(0,0)\cdot \frac{\partial \rho_M}{\partial u}\right|_o
= \alpha (0,0)  .                                      \tag4.7
$$
Therefore, 
$$
\left.\frac{\partial\delta}{\partial u}\right|_o = \alpha(0,0) > 0.
\tag4.8
$$
The Implicit Function Theorem now implies that the domain 
$F^\#(\Omega_M)$ admits a defining function $\rho_{N'}$
near $o$ given by
$$
\rho_{N'} (z,w) = \re w - G_{N'}(z,\bar z, v)           
$$
where 
$$
G_{N'}(z,\bar z, v) = \sum_{j,h\ge 1} b_{jh}(v) z^j \bar z^h
$$
is real analytic. 

Now, applying Hopf's Lemma,
we obtain that $\rho_{N'} \circ\tilde F$ is  a defining
function of $\Omega_{M'}$ near $o$, and more importantly that
it is in a \it reduced \rm form.  Thus, by letting
$$
\rho_{M'} \equiv \rho_{N'} \circ \tilde F,
$$
the assertion of the claim above is verified.  Consequently
the proof of 4.1 Proposition is complete.  \qed
\enddemo

Now we are ready to present the following main result of this section,
which shows the local rigidity of Siegel hypersurfaces
among the target hypersurfaces of proper holomorphic mappings from 
normalized weakly spherical hypersurfaces.  From now on, we denote by
$$
\varphi_k (z,w) = (z^k, w).
$$

\proclaim{4.3 Theorem} \it Let $(M,o)$ be a real analytic normalized
weakly spherical pointed CR hypersurface in $\C^2$ of order $k_0 > 1$.
Let $(\Sigma, o)$ be the pointed Siegel hypersurface given by
the defining equation
$$
\rho_{\Sigma} (z,w) \equiv \re w - |z|^2 = 0.
$$
If there is a holomorphic mapping $F:(M,o) \to (\Sigma, o)$
for which $o$ is a regular branch point, then
\roster
\item $(M,o)$ is defined by the equation $\re w - |z|^{2k_0} = 0$, and
\item $F(z,w) = (g_\Sigma\circ \varphi_{k_0}\circ g_M)(z, w)$, 
for some $g_\Sigma \in \aut (\Sigma, o)$ and $g_M \in \aut(M,o)$.
\endroster
\endproclaim

\demo{Proof}
By the preceding proposition, we may choose two elements 
$g_\Sigma$ $\in \aut (\Sigma, o)$ and $g_M$ $\in \aut(M,o)$ such that 
$$
\rho_M = \rho_\Sigma \circ g_\Sigma \circ F \circ g_M
$$
where $\rho_M$ is the weakly spherical normal form of $(M,o)$.
We put
$$
(g_\Sigma \circ F \circ g_M) (z,w) = f(z,w) = (f_1 (z,w), f_2 (z,w)).
$$
Then we have
$$
\re f_2 (z,w) - |f_1 (z,w)|^2
=
\re w - |z|^{2 k_0} 
- \sum_{j, h \ge 1 \atop j+h \ge 2 k_0 + 1} a_{jh}(\im w) z^j \bar z^h
                                                                                                                 \tag4.9
$$
where the coefficients $a_{jh}$ satisfy the vanishing conditions 
specified
in (3.4.2). 

Since $f(0,0) = (0,0)$, the Taylor expansions of $f_1$ and $f_2$ take 
the form
$$
\align
f_1 (z,w) &= b_{10} z + b_{01} w + R_{1,2} (z,w) \\
f_2 (z,w) &= \beta_{10} z + \beta_{01} w + R_{2,2} (z,w)
\endalign
$$
where $R_{k,m} (z,w)$ denotes the remainder term of order $m$ or 
higher in the Taylor expansion of $f_k (z,w)$.  Thus,
(4.9) becomes
$$
\multline
\re (\beta_{10} z + \beta_{01} w) +
\re R_{2,2}(z,w) - |b_{10} z + b_{01} w + R_{1,2}(z,w)|^2 \\
= \re w - |z|^{2k_0} - \sum a_{jh} (\im w) ~z^j \bar z^h    
\endmultline  \tag4.10
$$
Let $w=0$ in (4.10), then we get $\beta_{10} = 0$.  Then put $z=0$
in turn to get $\beta_{01} = 1$.  Since $k_0 > 1$, we must have
$\det Jf(0,0) = 0$.  This then implies $b_{10}=0$.  In
summary, we have
$$
b_{10} = \beta_{10} = 0, \quad \beta_{01} = 1           
$$
and 
$$
\re (R_{2,2} (z,w)) - |b_{01} w + R_{1,2} (z,w)|^2      
  = -  |z|^{2k_0} - \sum a_{jh}(\im w) z^j \bar z^h   \tag4.11
$$
Letting $z=0$ in (4.11), we obtain
$$
\re R_{2,2} (0,w) = |f_1 (0,w)|^2 .
$$
Since the left hand side is a harmonic function whereas the right 
hand side
attains its minimum at $w=0$, we deduce that
$$
\re R_{2,2}(0,w) = f_1 (0,w) = 0,
$$
for all values of $w$ in a neighborhood of $w=0$.  Thus,
in their Taylor expansions, neither $f_1$ nor $f_2$ has any nonzero
pure terms in $w^h$ with $h > 1$.  Moreover, it holds that
$$
b_{01} = 0. 
$$
By setting now $w=0$,
we get
$$ 
\multline
\re f_2 (z,0) - |b_{20} z^2 + \ldots + b_{k_0\,0} z^{k_0} + 
R_{1, k_0+1} (z,0)|^2 \\
= -|z|^{2k_0} - \sum a_{jh}(0) z^j \bar z^h            
\endmultline                                        \tag4.12
$$
where $R_{1, k_0+1} (z,0) = O(|z|^{k_0 + 1})$.

Comparing the terms in the power series expansions of both sides, we 
get
$$
b_{20} = \ldots = b_{(k_0-1)\, 0} = 0, 
\quad |b_{k_0\, 0}| = 1.    \tag4.13
$$
Moreover, we obtain $f_2 (z,0) \equiv 0$ since the monomials in the 
expansion of $\re f_2 (z,0)$ are types of either $z^\ell$ or $\bar z^\ell$
which are not found from any other terms of (4.12).  Therefore we now
have
$$
\multline
b_{k_0 \, 0} z^{k_0}\, \overline{R_{1,k_0+1} (z,0)} 
+ \overline{b_{k_0\,0} z^{k_0}}\, R_{1,k_0+1}(z,0) +
|R_{1,k_0+1}(z,0)|^2 \\
= \sum_{j,h \ge 1 \atop j+h \ge 2k_0 +1} a_{jh}(0) z^j \bar z^h 
\endmultline                                                   
\tag4.14
$$
where $|b_{k_0\, 0}| = 1$.

Now we use the vanishing condition (3.4.2) for $a_{jh}$ that is
$$
a_{k_0\,h} = 0 = a_{j\,k_0}, ~~\forall j, h \ge k_0 + 1.
$$
Applying this to (4.14) we obtain
$$
R_{1, k_0+1} (z,0) \equiv 0.
$$
As a consequence, we also have
$$
b_{j0} = 0, ~\forall j \ge k_0 + 1      \tag4.15
$$
and 
$$
a_{jh}(0) = 0, ~\forall j, h \ge 1 \text{ with } j+h \ge 2 k_0 + 1.  
\tag4.16
$$

Thus, it follows that
$$
\align
f_1 (z,w) & = b_{k_0\,0}\, z^{k_0} + r_1 (z,w) \\
f_2 (z,w) & = w + r_2 (z,w)
\endalign
$$
where both $r_1(z,w)$ and $r_2(z,w)$ are holomorphic functions 
of the class $O(|zw|)$.

Now, we let $w = \epsilon \in \R$, with $\epsilon$ sufficiently small
but positive.  Using (4.16), the identity (4.10) now becomes
$$
\epsilon + \re r_2 (z,\epsilon) - |b_{k_0\,0} z^{k_0} + r_1
(z,\epsilon)|^2
= \epsilon - |z|^{2k_0}
\tag4.17a
$$
or, equivalently
$$
\re r_2(z,\epsilon) - b_{k_0\,0} z^{k_0}\, \overline{r_1 (z,
\epsilon)}
- \overline{b_{k_0\,0} z^{k_0}}\, r_1 (z,\epsilon) = |r_1
(z,\epsilon)|^2
\tag4.17b
$$
Here, we want to show first that $r_2 (z,\epsilon) \equiv 0$. For
each fixed $\epsilon > 0$ that is sufficiently small, we expand each
term in the identity (4.17b) using the Taylor expansion of $r_1$ and
$r_2$ in the variable $z$ at the origin.  The monomial terms in the 
expansion of $\re r_2(z,\epsilon)$ consist solely of the form $z^\ell$
and $\bar z^\ell$, whereas the remaining terms contain only the mixed
power terms of $z$ and $\bar z$.  Thus the comparison of coefficients
immediately implies that $r_2 (z, \epsilon) = 0$ for every $z$ in an
open neighborhood of the origin and for every sufficiently small
$\epsilon > 0$.  Since $r_2$ is holomorphic in $z,w$, this implies
immediately that $r_2(z,w)$ is indeed identically zero.

Due to the preceding arguments, (4.17a) now becomes
$$
|b_{k_0\,0} z^{k_0} + r_1(z,\epsilon)|^2 = |z|^{2k_0}
$$
or equivalently
$$
\left| b_{k_0\,0} + \frac{r_1(z,\epsilon)}{z^{k_0}} \right| = 1
$$
for all $z$ in a neighborhood of the origin in $\C$ and for all
sufficiently small values of $\epsilon > 0$.
It then follows immediately that
$$
b_{k_0,0} z^{k_0} + r_1(z,\epsilon) = e^{ig(\epsilon)} z^{k_0}
$$
for some real-valued real-analytic function $g$ defined
in an open neighborhood of the origin in $\R$. Thus, we conclude that 
$r_1(z,w) = z^{k_0}\hat r(w)$, where $\hat r(w)$ is holomorphic in
the complex variable $w$. 

Therefore, we now arrive at that our mapping $f$ is of the form 
$$
f(z,w) = (z^{k_0}(b_{k_0,0}+\hat r(w)) , w).
$$
We apply this information to the identity (4.9) and  
we immediately observe the following: 
\roster
\item all coefficients $a_{ij}$ vanish identically, and
\item $|b_{k_0,0}+\hat r(w)| \equiv 1$, and hence 
 $b_{k_0,0}+\hat r(w) \equiv \lambda$, for some constant $\lambda \in \C$
with $|\lambda| = 1$.
\endroster
Thus, we obtain 
$$
(g_\Sigma \circ F \circ g_M) (z,w) = f(z,w) = (\lambda z^{k_0}, w),
$$
which proves the assertion.
\qed
\enddemo

\head
5. Factorization of proper maps through a complex ellipsoid
\endhead

Denote by
$$
E_m \equiv \{(z,w) \in \C^2 \st |z|^{2m} + |w|^2 = 1 \}
$$
and by
$$
\Omega_{E_m} \equiv \{(z,w) \in \C^2 \st |z|^{2m} + |w|^2 < 1 \}
$$
the \it complex ellipsoid \rm and the \it Thullen domain of order 
$m$\rm, respectively.  Similarly, denote by
$$
\Sigma_m \equiv \{(z,w) \in \C^2 \st \re w - |z|^{2m} = 0 \}
$$
and 
$$
\Omega_{\Sigma_m} \equiv \{(z,w) \in \C^2 \st \re w - |z|^{2m} > 0 \}
$$
the \it normalized Siegel hypersurface \rm and the \it Siegel upper-half space
\rm of order $m$.
Recall that, for any positive integer $m$, the map
$$
\mu_m (z,w) = \left( \frac{z}{(1+w)^{1/m}}, \frac{1-w}{1+w} \right)
$$
sends {\it biholomorphically\/} the pointed CR hypersurface
$(E_m,(0,1))$ 
to the pointed CR hypersurface $(\Sigma_m,(0,0))$. 

For simplicity, in the following we will
often write just $E = \partial B^2$ instead of $E_1$,
 $\Sigma$ instead of $\Sigma_1$ and $\mu$ instead of $\mu_{1}$.
We also continue using the notation
$\varphi_m(z,w) = (z^m,w)$.  Observe that, for any positive 
integer $m$
$$
\varphi_m \circ \mu_m = \mu\circ \varphi_m    \tag5.1
$$

We now present a key lemma.

\proclaim{5.1 Lemma} \it
Let $U_1$ and $V_1$ be open neighborhoods of the point $(0,1)$
in $\C^2$ and let $F:U_1 \to V_1$ be a local biholomorphism
of the pointed CR hypersurface $(E_m, (0,1))$ into itself. Then $F$
extends uniquely to a global automorphism of the Thullen
domain $\Omega_{E_m}$.         
\endproclaim

\demo{Proof} Let $\mu_m$ and $\varphi_m$ be as above.
Furthermore, let us denote by
$$
\hat F = \mu_m \circ F\circ \mu_m^{-1} : \Sigma_m \to \Sigma_m
$$ 
and 
$$
G = \varphi_m \circ \hat F : \Sigma_m \to \Sigma
$$
An important step toward the proof is the following

\roster
\item"{\bf Claim.}" \it There exists a local
automorphism $\psi$ of $(\Sigma,o)$ such that
$$
\rho_{\Sigma_m} = \rho_\Sigma \circ \psi \circ G
$$
where $\rho_{\Sigma_m}$ and $\rho_\Sigma$ are the defining functions
in normal form for $\Sigma_m$ and $\Sigma$ respectively.
\rm
\endroster

\it Proof of the claim: \rm Notice that $\hat F$ preserves the set 
$\{ z = 0\} \cap (\Sigma_m, (0,0))$, 
which is the set of weakly pseudoconvex points of $\Sigma_m$. 
Hence, it follows that
$$
\hat F(z,w) = (z k(z,w), f(w) + z g(z,w))
$$
for some holomorphic functions $k$, $f$ and $g$.  We rewrite $g(z,w)$ in
its MacLaurin expansion:
$$
g(z,w) = \sum g_{r,s} z^r w^s
$$
We want to show that $g\equiv 0$ first.  Let us restrict $\hat F$ 
to the set of points of the type $(z,|z|^{2m} + i c)$, for sufficiently
small positive real values of $c$.  Since $\hat F$ is a local automorphism
of $\Sigma_m$ preserving $o=(0,0)$, we have
$$
\re(z g(z, |z|^{2m} + ic))
= |z|^{2m} | k(z, |z|^{2m} + ic) |^{2m} -
\re f(|z|^{2m} + ic)                          \tag5.2
$$
Now, we apply the operator 
$\left.\frac{\partial^{r+1}} {\partial z^{r+1}}\right|_{z=0}$
to both sides of (5.2).  Then we deduce that there exists $\eta > 0$
such that 
$$
\sum_s (r+1)! (ic)^s g_{rs} \equiv 0
$$
for all $0 < c <\eta$.  This implies $g_{rs} = 0$ for all $r, s$.
Hence, $g \equiv 0$.

Now we have 
$$
\hat F(z,w) = (z k(z,w), f(w))
$$
Rewriting the identity (5.2) with $g$ identically zero, we get
$$
\re f(|z|^{2m}) =
|z|^{2m} | k(z, |z|^{2m})|^{2m} \quad, \quad \forall z       \tag5.3
$$
Since $\hat F(0,0)=(0,0)$, it follows that $f(0)=0$.  Hence, when we rewrite
$f$ in its MacLaurin expansion, we get
$$
f(w) = \sum_{j \geq 1} f_j w^j
$$
It follows by (5.3) that 
$$
0 \equiv \frac{1}{2} \sum_{j\geq 1} (f_j + \overline f_j)
|z|^{2m(j -1)} - |k(z,|z|^{2m})|^{2m}                       \tag5.4
$$
From this identity, we may deduce in particular that
$$
0 =
\left.\frac{\partial^h}{\partial z^h}
|k(z,|z|^{2m})|^{2m}\right|_{z=0}
$$
for any positive integer $h$.  Notice also that
$$
\frac{\partial^h}{\partial z^h}
|k(z,|z|^{2m})|^{2m} =
\left.\frac{\partial^h}{\partial t^h}
|k(t,|z|^{2m})|^{2m}\right|_{t=z} + O(|z|)                \tag5.5
$$
Since $\hat F$ is a local automorphism at $(0,0)$, $k(0,0)\neq 0$.
Altogether, it follows by (5.5) that
$\frac{\partial^h}{\partial t^h}k(t,0) = 0$ 
for any positive integer $h$.
This implies that $k(z,w) = k(w)$.  In conclusion, the map
$\hat F$ is of the form
$$
\hat F(z,w) = (z k(w), f(w))
$$
with $k(0) \neq 0$ and $\frac{\partial f}{\partial w}(0) \neq 0$.
Thus, 
$$
G(z,w) = (z^m k(w)^m, f(w))
$$
Let us consider the local 
biholomorphism $\psi$ in a neighborhood of the origin in $\C^2$
defined by
$$
\psi^{-1}(z,w) = (z k^m(w), f(w))
$$
Also denote by
$$
G^\#(z,w) = \psi\circ G(z,w)
$$
Then
$$
G^\#(z,w) = (z^m, w) = \varphi_m (z,w)    \tag5.6
$$
Notice that 
$
G^\# ( \Sigma_m \cap U_{\Sigma_m}) = \varphi_m (\Sigma_m \cap 
U_{\Sigma_m}) \subset \Sigma
$
and that $\rho_{\Sigma}\circ G^\# = \rho_\Sigma \circ \varphi_m
= \rho_{\Sigma_m}$.  Moreover, it is easy to see that $\psi$ is 
a local automorphism of $\Sigma$, because the images of $G$ and 
$G^\# = \psi\circ G$ are both contained in $\Sigma$.  

Consequently, we arrive at
$$
\rho_{\Sigma_m} = \rho_{\Sigma}\circ \psi\circ G  
$$
and the claim is proved.
\medskip

Now we want to complete the proof of Lemma 5.1.

By \cite{1}, $\psi$ extends uniquely to a global 
automorphism, say $\Psi$, of $\Sigma$ fixing $o$. 
By (5.6) we have
$$
\Psi\circ G = \varphi_m
$$
This, together with (5.1), leads us to the identity
$$
\varphi_m \circ F = \mu^{-1}\circ \Psi^{-1}
\circ \mu\circ \varphi_m = \Phi\circ \varphi_m          \tag5.7
$$
where $\Phi = \mu^{-1}\circ \Psi^{-1} \circ \mu$ is a global 
automorphism of $B^2$.  Moreover, 
$$
\Phi (\{(0,it) \st t \in \R\}) \subset \{(0,it) \st t \in \R\}.
$$
Appealing to an explicit formula for the automorphism $\Phi$, we are
immediately able to deduce that each branch of 
$\varphi_m^{-1} \circ \Phi \circ \varphi_m$ defines a biholomorphism
of whole $\Sigma_m$.  So the proof of the lemma is complete by letting 
$F$ be one of the branches.
\qed
\enddemo

Now we present the main theorem of this paper.  

\proclaim{5.2 Theorem} \it
Let $D$ be a bounded simply connected pseudoconvex domain in $\C^2$ 
with a real analytic boundary.  Assume that $f:D \to B^2$ is a 
proper holomorphic mapping with generic degree $m$ 
such that the analytic variety $Z_{df}$ admits an 
irreducible component $V$ satisfying: 
\roster
\item $f^{-1}(f(V \cap \partial D)) = V \cap \partial D$;
\item $V \cap \partial D$ is connected and contains no singular
point of the variety $V$.
\endroster
Then, there exists a proper holomorphic mapping $g:D \to \Omega_{E_k}$ 
from $D$ to the Thullen domain $\Omega_{E_k}$,
where $k$ is the generic degree of $f$ in a sufficiently small
tubular neighborhood of $V\cap \partial D$, 
which extends holomorphically across $\partial D$ such that
$$
f = \beta \circ \varphi_k \circ g
$$
for some holomorphic automorphism $\beta$ of $B^2$.
\endproclaim

Notice that this immediately implies

\proclaim{5.3 Theorem} \it
\it Let $D$ be a bounded simply connected pseudoconvex domain in $\C^2$
with a real analytic boundary.  If $D$ admits a generically
$k$-to-1 proper holomorphic mapping $f:D \to B^2$ such that there  
exists an irreducible subvariety $V$ of $Z_{df}$ satisfying:
\roster
\item $f$ is a local $k$-to-1 branched covering at every point of 
$V \cap \partial D$;
\item $V \cap \partial D$ is connected and contains no singular
point of the variety $Z_{df}$;
\endroster
then $D$ is biholomorphic to $\Omega_{E_k}$.
\endproclaim

The rest of the section is devoted to the proof of 5.2 Theorem.  We 
first consider the concept of holomorphic correspondences.

A \it holomorphic correspondence \rm from a domain $\Omega$ in 
$\C^m$ to $\Omega'$ in $\C^n$ is a complex analytic set 
$S$ in $\Omega \times \Omega'$ such that $\pi_\Omega(S)=\Omega$,
where $\pi_\Omega : \Omega \times \Omega' \to \Omega$ is the 
standard projection onto $\Omega$.  We denote the correspondence by
$$
S: \Omega \multimap \Omega'
$$

For a holomorphic correspondence $S:\Omega \multimap \Omega'$,
we denote by  
$$
S^{-1} (G) = \pi_\Omega (\pi_{\Omega'}^{-1} (G) \cap S), 
  \text{ for } G \subset \Omega'
$$
Furthermore, we call $S$ \it proper \rm if $S^{-1}(K)$ is compact for
every compact subset $K$ of $\Omega'$.

An important concept associated with the holomorphic correspondences
concerns whether it is actually realized by a union of the graphs of
holomorphic mappings.  Precisely speaking, a proper holomorphic 
correspondence $G : \Omega \multimap \Omega'$ between bounded domains
in $\C^n$, which extends holomorphically across $\partial\Omega$ is 
said
to \it split locally at $p_0 \in \overline{\Omega}$\rm, if there exist
a neighborhood $U_0$ of $p_0$ in $\C^n$ and $m$ holomorphic maps
$f_j : U_0 \to f_j(U_0), ~(j=1,\ldots, m)$ such that
$$
G|_{U_0} 
= \bigcup_{j=1}^m \{ (z,f_j(z)) \st z \in U_0 \} 
$$
A holomorphic correspondence is said to \it split globally\rm, if it is
the union of the graphs of globally defined holomorphic mappings.  In
our proof here, we use the following fact that was observed in
Lemmata 3.6 and 3.7 by Bedford-Bell (\cite{5}):

\it Let $\Omega$ be a bounded simply connected domain, and
let $\Omega'$ be a bounded domain in $\C^n$.  
Let $G : \Omega \multimap \Omega'$ be a proper holomorphic 
correspondence
that extends holomorphically across the boundary of $\Omega$.  Then
$G$ splits globally if and only if it splits locally at every boundary
point of $\Omega$. \rm

Now, let $f:D \to B^2$ be a proper holomorphic mapping given in the
hypothesis of 5.2 Theorem.

By the extension theorem of Diederich-Forn{\ae}ss \cite{9}, the mapping 
$f$ extends holomorphically to an open neighborhood of the closure 
$\overline{D}$ of $D$ onto a neighborhood of $B^2$ mapping $\partial D$ onto
$\partial B^2$.  We consider the holomorphic correspondence
$$
G \equiv \alpha^{-1}\circ \varphi_k^{-1} \circ \beta \circ u_p\circ f : 
\overline{D} \multimap \overline{\Omega}_{E_k}
$$
for some appropriate $\alpha \in \aut {\Omega}_{E_k}$ 
and $\beta, u_p \in \aut B^2$, which are to be 
determined later. In order to verify the assertion
of 5.2 Theorem, it is enough to show that the proper holomorphic
correspondence $G$ above splits at every boundary point of $D$.
Then, $G$ will consist of $k$ graphs of holomorphic mappings, each
of which will provide the desired factorization map.

We will make an appropriate choice of 
$\alpha \in \aut {\Omega}_{E_k}$ 
and $\beta, u_p\in \aut B^2$  
as well as check the local splitting of $G$ at every point of 
$V \cap \partial D$.  Then we will check the local 
splitting of $G$ at the points not in $V \cap \partial D$.

\subhead
\it
Step 1.
Choice of $\alpha$, $\beta$ and $u_p$  \rm
\endsubhead

Pick a point $p \in V \cap \partial D$ and choose a unitary map 
$u_p \in \aut B^2$
so that $u_p \circ f(p) = (0,1) \in \partial B^2$.  Then, using 
the notation introduced at the beginning of this 
section,  consider the
biholomorphic mapping $\mu:B^2 \to \Omega_{\Sigma}$.
Now, choose a sufficiently small open neighborhood $U_p$ of $p$ in 
$\C^2$
so that the normalization by Barletta-Bedford \cite{4} can be 
applied.  Namely, the boundary $\partial D$ in $U_p$ is weakly
spherical and there exists a local biholomorphic mapping $\psi_p$
from $U_p$ onto an open neighborhood of the origin such that
$\psi_{p}(\partial D)$ is now represented by the normal form
$$
\re w = |z|^{2k} + \text{higher order terms}
$$
where the higher order terms satisfy the conditions specified in
(3.4.1) and (3.4.2).  
Composing these maps, we arrive at the holomorphic mapping
$$
F \equiv \mu \circ u_p \circ f \circ \psi_p^{-1} :
(\psi_p (\partial D), o) \to (\Sigma, o)
$$
from a normalized weakly spherical pointed CR surface to the 
normalized Siegel pointed CR surface.  

By 4.3 Theorem, we deduce that
$\psi_p (\partial D\cap U_p)$ is in fact a neighborhood of the origin 
in the hypersurface
$\Sigma_k$
and there exist $\hat \beta \in \aut (\Sigma,o)$ 
and $\hat \alpha \in \aut(\Sigma_k,o)$ such that 
$$
F(z,w) = (\hat\beta^{-1} \circ \varphi_k \circ 
\hat\alpha)(z,w) .                                     \tag 5.8
$$
By [1] and 5.1 Lemma the mappings $\mu^{-1}\circ\hat\beta\circ\mu$ 
and $\mu_k^{-1}\circ \hat\alpha\circ\mu_k$ extend, 
respectively, to $\beta \in Aut(B^2)$ and $\alpha\in Aut(\Omega_{E_k})$.
This implies, in particular, that $\hat\alpha$
and $\hat\beta$ also extend to global automorphisms
of $\Omega_{\Sigma}$ and $\Omega_{\Sigma_k}$, respectively.

Define a global holomorphic correspondence
$$
G = \alpha^{-1}
\circ\varphi^{-1}_k\circ \beta
\circ u_p \circ f                                     \tag5.9
$$
Note that in $U_p$, $G$ has a local expression
$$
G=
\alpha^{-1} \circ \varphi_k^{-1} \circ
\beta \circ \mu^{-1} \circ F\circ \psi_p               
$$
Using (5.8) and the definition of $\beta$, we get
$$
G = \alpha^{-1}\circ \varphi_k^{-1} 
\circ\mu^{-1}\circ\varphi_k\circ\hat\alpha \circ\psi_p
$$
It follows by (5.1) that 
$$
G = \alpha^{-1} \circ \varphi_k^{-1} \circ \varphi_k \circ \alpha
\circ \mu_k^{-1} \circ \psi_p                                \tag5.10
$$
Since the correspondence $\varphi^{-1}_k\circ \varphi_k$ splits, 
and since $\alpha$, $\mu_k$ and $\psi_p$ are local biholomorphisms,
it follows by (5.10) that the correspondence $G$ splits in a
small neighborhood of $p$ where $\psi_p$ is defined.  However, we
have yet to see if $G$ splits at every boundary point.  We will 
see this in the next two and final steps.

\subhead
\it Step 2. Splitting of $G$ at points of $\partial D\cap V$ \rm
\endsubhead

Let $q \in \partial D\cap V$.  Assume for a moment that the 
neighborhood $U_q$ of $q$ satisfies the following additional condition
$$
U_q\cap U_p\neq \emptyset.
$$
In $U_p\cap U_q$, (5.10) is valid and we have
$$
G = \alpha^{-1} \circ \varphi_k^{-1}\circ
\varphi_k \circ \alpha \circ \hat\psi_p \circ \hat\psi_q^{-1}
\circ \hat\psi_q                                        \tag5.11
$$
where $\hat\psi_p = \mu_k^{-1} \circ \psi_p$ and
$\hat\psi_q = \mu_k^{-1} \circ \psi_q$, and where $\psi_q$ is 
the Barletta-Bedford normalization map for $\partial D \cap U_q$.
Since $\hat\psi_p\circ \hat\psi^{-1}_q$ extends
to an automorphism of $\Omega_{E_k}$ by 5.1 Lemma, the expression
(5.11) is well-defined on $U_q$.  This shows that $G$ splits
at every point of $U_p \cup U_q$.

For an arbitrary point of $\partial D \cap V$, an inductive repetition 
of this argument yields the desired conclusion for the current step, 
because $\partial D\cap V$ is a compact connected set.

\subhead
\it Step 3. Splitting of $G$ at points of $\partial D\setminus V$ \rm
\endsubhead

Let $q \in \partial D\setminus V$. Recall 
$$
G = \alpha^{-1}\circ \varphi_k^{-1}\circ \tilde f
$$
where $\alpha$ is a biholomorphism of the Thullen domain
$\Omega_{E_k}$ and $\tilde f = \beta \circ u_p \circ f$ is
a proper holomorphic map branching at every point of $V$.
Thus, for $q \in \partial D\setminus V$, we have that $\tilde f(q)$ 
is a strictly pseudoconvex point that does not belong to
$\tilde f(V\cap \partial D) \subset \{ z = 0\}$, by (1) in the 
hypothesis of the theorem.  Thus, $G(q)$ consists only of
strictly pseudoconvex points. The splitting of $G$
is then guaranteed by the proof of Theorem 3 of [5].  
\qed

\enddocument